\documentclass[11pt]{article}
\usepackage[english]{babel}
\usepackage{amsmath,amsthm}
\usepackage{amsfonts}
\usepackage[comma,numbers,square,sort&compress]{natbib}
\usepackage{graphicx,subfigure,amsmath,amssymb,mathrsfs,amsfonts,amstext,amsthm}
\usepackage{latexsym, amssymb}
\usepackage{graphicx}
\usepackage{subfigure}
\usepackage{pgf,tikz}
\usepackage{pstricks}
\newtheorem{thm}{Theorem}[section]

\newtheorem{lem}[thm]{Lemma}

\theoremstyle{definition}
\newtheorem{defn}[thm]{Definition}
\theoremstyle{remark}

\numberwithin{equation}{section}
\begin{document}

\title{\bfseries\textrm{Robust Fragmentation Modeling of Hegselmann-Krause-Type Dynamics*}
\footnotetext{*This research was supported by the National Natural Science Foundation of China under grants No. 61671054, 61773054, 91427304, 61333001, and 11688101, and by the National Key Basic Research Program of China (973 program)
under grant 2014CB845301/2/3 and the Fundamental Research Funds for the Central Universities under grant No. FRF-TP-17-087A1.\\
\,\,W. Su, J. Guo, X. Chen are with School of Automation and Electrical Engineering, University of Science and Technology Beijing \& Key Laboratory of Knowledge Automation for Industrial Processes, Ministry of Education, Beijing 100083, China, {\tt suwei@amss.ac.cn, guojin@amss.ac.cn, cxz@ustb.edu.cn}. G. Chen is with National Center for Mathematics and Interdisciplinary Sciences \& Key Laboratory of Systems and
Control, Academy of Mathematics and Systems Science, Chinese Academy of Sciences, Beijing 100190,
China, {\tt chenge@amss.ac.cn}.}}
\author{Wei Su \and Jin Guo \and Xianzhong Chen \and Ge Chen}

%
\date{}%
\maketitle
\begin{abstract}
In opinion dynamics, how to model the enduring fragmentation phenomenon (disagreement, cleavage, and polarization) of social opinions has long possessed a central position. It is widely known that the confidence-based opinion dynamics provide an acceptant mechanism to produce fragmentation phenomenon. In this study, taking the famous confidence-based Hegselmann-Krause (HK) model, we examine the robustness of the fragmentation coming from HK dynamics and its variations with prejudiced and stubborn agents against random noise. Prior to possible insightful explanations, the theoretical results in this paper explicitly reveal that the well-appearing fragmentation of HK dynamics and its homogeneous variations finally vanishes in the presence of arbitrarily tiny noise, while only the HK model with heterogenous prejudices displays a solid cleavage in noisy environment.
\end{abstract}


\section{Introduction}
Opinion dynamics has recently displayed its increasing attraction to the researchers from diverse areas, including sociology, mathematics, information science, physics and so on \cite{Castellano2009,Proskurnikov2017}. One quantitative way to reveal the opinion behaviors is to technically model its evolution mechanism, and various agent-based models have been established during the past decades \cite{DeGroot1974,Friedkin1999,Deffuant2000,Krause2000,Hegselmann2002,Friedkin2016,Zhang2013,Zhang2014}. Roughly categorized into two classes, namely topology-based and confidence-based, all these models succeeded in capturing some features of opinion evolution, and demonstrated plenty of behavioral pictures coinciding with social reality, ranging from the basic agreement or disagreement of the opinions \cite{DeGroot1974,Hegselmann2002}, to the more complex autocratic and democratic social structures \cite{Jia2015}.

One of predominantly important issues concerned in social dynamics is how to build the model to produce the social disagreement or community cleavage \cite{Hegselmann2002,Abelson1964,Friedkin2015,Kurahashi2016}. In \cite{Krause2000,Hegselmann2002}, a notable model widely known as Hegselmann-Krause (HK) opinion dynamics later, was developed based on bounded confidence and presented well-formed fragmentation. Meanwhile, G. Deffuant etc. established a similar-spirit model along the same line \cite{Deffuant2000}.
Bounded confidence, which limits the neighbor opinions that agents take into account when updating their opinions, ever since has been believed to be markedly capable of generating fragmentation of opinion dynamics. Afterwards the HK model got fully exploited and fruitful results were theoretically obtained \cite{Lorenz2005,Blondel2009}. Other than bounded confidence, prejudiced and stubborn agents are also found as two remarkable factors that take group cleavage as outcomes \cite{Hegselmann2006,Friedkin2015,Dabarera2016,Proskurnikov2017prejudice}.

However, whether bounded confidence can produce robust clusters has been challenged a lot (see \cite{Kurahashi2016} or references therein), and a very late study theoretically reveals that the opinions under HK dynamics will spontaneously achieve a consensus in the presence of random noise, no matter whatever initial opinions and tiny noises were given \cite{Su2015}. In HK model, each agent updates its states by averaging the opinions of its confidence-dependent neighbors. This local rule of self-organization in HK model allows the fragmentation behavior to emerge. When random noise is admitted, the interaction of opinions occurs in a stochastic yet more connected way, and in consequence leads to a consensus. Moreover, in a variation of HK model with actually the homogeneously prejudiced agents, the same phenomenon was verified, where driven by random noise, the ever divergent opinions get merged to the prejudice value \cite{Su2017free}. With these affirmative facts of noise-induced consensus in hand, we can conclude, as a byproduct of these facts, that the elegant HK dynamics alone or with homogeneously prejudiced agents is inadequate to produce the fragmentation phenomenon in noisy environment. However, random noise arises ubiquitously in natural and social systems \cite{Sagues2007,Mas2010,Pineda2011,Guo2017}, whence how the elegantly structured HK dynamics brings about robust fragmentation against random noise requires further exploration.

In this paper, we propose and then examine the noisy HK models with prejudiced and with stubborn agents in both homogeneous and heterogeneous cases. Just as the HK model with homogeneously prejudiced agents uncovered in \cite{Su2017free}, it will be rigorously proved that the HK model with homogeneously stubborn agents fails likewise to exhibit a robust fragmentation against random noise, with the noisy opinions finally getting synchronized to the stubborn agents. While the noisy HK model with heterogeneously stubborn agents will only partly show community cleavage for some specific initial opinions, it is finally established that the HK dynamics with heterogeneously prejudiced agents will demonstrate a markedly robust fragmentation in the presence of noise. This discovery confirms a plainly intuitive idea that \emph{\textbf{it is the innate difference lying within the community rather than the mere bounded confidence of individuals at that accounts for the ubiquitous social cleavage.}}

The rest of the paper is organized as follows: Sec. \ref{Secmodel} introduces noisy HK models and their spontaneous consensus in the presence of noise; Sec. \ref{Sechetermodel} presents the theoretical discoveries of the noisy HK model with heterogeneous prejudices; Sec. \ref{Secsimul} provides some numerical simulations to verify the theoretical conclusions, and finally Sec. \ref{Secconc} concludes the paper.

\section{$\phi$-consensus of HK-type models}\label{Secmodel}
\renewcommand{\thesection}{\arabic{section}}

In this part, we will examine the fragmentation of noisy HK model and its variations with homogeneously prejudiced and stubborn agents. The theoretical results will show that the well-appearing fragmentation of these models disappears in the presence of arbitrarily tiny noise.
To begin with, the noisy HK models and strict definition of consensus in the noisy case are need.
\subsection{Noisy HK models and $\phi$-consensus}
Suppose there are $n$ agents in the group with $\mathcal{V}=\{1,2,\ldots,n\}$, $x_i(t)\in[0,1], i\in \mathcal{V}, t\geq 0$ is the opinion value of agent $i$ at time $t$, then
the basic noisy HK model is adopted following that in \cite{Su2015}:
\begin{equation}\label{HKnoise0}
  x_i(t+1)=\left\{
           \begin{array}{ll}
             1,  & \hbox{~$x_i^*(t)>1~$} \\
             x_i^*(t),  & \hbox{~$x_i^*(t)\in[0, 1]~$} \\
             0,  & \hbox{~$x_i^*(t)<0~$}
           \end{array}~\forall i\in\mathcal{V},\, t\geq 0,
         \right.
\end{equation}
where
\begin{equation}\label{neigh0}
  \mathcal{N}(i,x(t))=\{j\in \mathcal{V}:|x_j(t)-x_i(t)|\leq \epsilon\}
\end{equation}
and
\begin{equation}\label{xit0}
  x_i^*(t)=|\mathcal{N}(i, x(t))|^{-1}\sum\limits_{j\in \mathcal{N}(i, x(t))}x_j(t)+\xi_i(t+1)
\end{equation}
with $\epsilon\in(0,1]$ being confidence bound, and $\{\xi_i(t)\}_{i\in\mathcal{V},t> 0}$ the random noises which are i.i.d. with $E\xi_1(1)=0, E\xi_1^2(1)>0$ and $|\xi_1(1)|\leq \delta$ a.s. for $\delta\geq 0$.

Without noise in (\ref{xit0}) (i.e. $\delta=0$), system (\ref{HKnoise0})-(\ref{xit0}) degenerates to the original HK model based on bounded confidence, which will converge in finite time in the sense that, there exist $T\geq 0, x_i^*\in[0,1], i\in\mathcal{V}$ such that $x_i(t)=x_i^*, t\geq T$. If $x_i^*=x_j^*$ for any $i,j\in\mathcal{V}$, the opinions is said to reach consensus; otherwise, there must exist $i,j\in\mathcal{V}$ with $|x_i^*-x_j^*|>\epsilon$ and fragmentation forms.

Other than bounded confidence, prejudiced and stubborn agents are found to be another two remarkable causes of well-formed fragmentation combined with HK dynamics. In this part, we will first explore the noisy HK models with homogeneously prejudiced and stubborn agents. Here, ``homogeneous'' refers to that all the prejudiced agents or the stubborn ones possess an identically constant opinion value.


The noisy HK model with homogeneously prejudiced agents is obtained by modifying (\ref{xit0}) as follows:
\begin{equation}\label{xit1}
\begin{split}
  x_i^*(t)=(1-\alpha I_{\{i\in\mathcal{S}_1\}})\sum\limits_{j\in\mathcal{N}_i(x(t))}\frac{x_j(t)}{|\mathcal{N}_i(x(t))|}+\alpha I_{\{i\in\mathcal{S}_1\}}J_1+\xi_i(t+1).
\end{split}
\end{equation}
where $\mathcal{S}_1\subset\mathcal{V}$ is the set of prejudiced agents, $J_1\in[0,1]$ is the prejudice value, and $\alpha\in(0,1]$ is the attraction strength of prejudice value. As usual, without noise the model generates conspicuous fragmentation \cite{Hegselmann2006,Su2017free}.

For the noisy HK model with homogeneously stubborn agents, to be convenience, we introduce additionally in the model (\ref{HKnoise0})-(\ref{xit0}) a stubborn agents set, $\it \mathcal{B}_1$,  whose opinion values satisfy
\begin{equation}\label{stubbornvalues1}
  x_i(t)\equiv B_1,\quad i\in\mathcal{B}_1, B_1\in[0,1]
\end{equation}
and the neighbor set in (\ref{neigh0}) is modified as
\begin{equation}\label{neigh1}
  \mathcal{N}(i,x(t))=\{j\in \mathcal{V}\cup\mathcal{B}_1:|x_j(t)-x_i(t)|\leq \epsilon\}.
\end{equation}
Then (\ref{HKnoise0}), (\ref{xit0})-(\ref{neigh1}) describe a noisy HK model with homogeneously stubborn agents. Likewise, when the noise strength $\delta=0$, the system degenerates to the noise-free one, which is claimed converging and also could produce fragmentation \cite{Chazelle2015}.


\subsection{$\phi$-consensus}
Due to persistent disturbance of noise, the definition of consensus of noisy HK models slightly differs:
\begin{defn}\label{robconsen}\cite{Su2016}
Define
\begin{equation*}\label{opindist}
  d_{\mathcal{V}}(t)=\max\limits_{i, j\in \mathcal{V}}|x_i(t)-x_j(t)|~~\mbox{and}~~d_{\mathcal{V}}=\limsup\limits_{t\rightarrow \infty}d_{\mathcal{V}}(t).
\end{equation*}
For $\phi\in[0,\infty)$,\\
(i) If $d_{\mathcal{V}} \leq \phi$, we say the system will reach $\phi$-consensus.\\
(ii) If $P\{d_{\mathcal{V}} \leq \phi\}=1$, we say almost surely (a.s.) the system will reach $\phi$-consensus.\\
(iii) Let $T=\inf\{t: d_{\mathcal{V}}(t')\leq \phi \mbox{ for all } t'\geq t\}$.
 If $P\{T<\infty\}=1$, we say a.s. the system reaches $\phi${-\it consensus} in finite time.

 Furthermore for a constant $A\in[0,1]$, define $d_\mathcal{V}^{A}(t)=\max\limits_{i\in\mathcal{V}}|x_i(t)-A|$ and $d_\mathcal{V}^A=\limsup\limits_{t\rightarrow\infty}d_\mathcal{V}^{A}(t)$, if substitute $d_\mathcal{V}^A$ in (i)-(iii) for $d_{\mathcal{V}}$, the system is accordingly said to reach $\phi$-{\it consensus with} $A$.
\end{defn}
%


\subsection{Spontaneous $\phi$-consensus caused by noise}
Though the HK model with either prejudiced or stubborn agents enables a well-appearing fragmentation phenomenon, the introduction of random noise, even very tiny, makes the $\phi$-consensus of ever-divergent opinions spontaneously emerge. To be specific, we have
\begin{thm}\label{thmhomo}
Given any $x(0)\in[0,1]^n, \epsilon\in(0,1]$:\\
(a) for all $\delta\in(0,\epsilon/2]$, a.s. the system (\ref{HKnoise0})-(\ref{xit0}) will reach $2\delta$-consensus in finite time.\\
(b) for all $\delta\in(0,\underline{\delta}]$, a.s. the system (\ref{HKnoise0}), (\ref{neigh0}), (\ref{xit1}) will reach $\bar{\delta}$-consensus with $J_1$. Here, $\bar{\delta}$ and $\underline{\delta}$ are constants determined by parameters $n,\epsilon,\alpha, |\mathcal{S}_1|$ (see \cite{Su2017free}).\\
(c) for all $\delta\in(0,\frac{\epsilon}{2(n+1)})$, a.s. the system (\ref{HKnoise0}), (\ref{xit0})-(\ref{neigh1}) will reach $2\delta$-consensus and meanwhile $(n+1)\delta$-consensus with $B_1$;\\
\end{thm}

Theorem \ref{thmhomo} evidently shows that in the presence of tiny noise ($\delta$ is teeny), the HK-based opinion dynamics with homogeneous variations lose their ability of generating fragmentation phenomenon, and almost surely achieve a fairly approximate consensus (see Fig. \ref{stubbornsfig}). We did not investigate with theoretical analysis the opinion behaviors when the noise strength is large, where even the gathered opinions could get divergent and tanglesome, as verified in Theorem 8 of ~\cite{Su2015}. Conclusions (a) and (b) follow directly from Theorem 2 of ~\cite{Su2015} and Theorem 3.1 of ~\cite{Su2017free}. Now we proceed to the
the proof of conclusion (c). To begin with, some preliminary lemmas are need.
\begin{lem}\cite{Blondel2009}\label{monosmlem}
Suppose $\{z_i, \, i=1, 2, \ldots\}$ is a nondecreasing (nonincreasing) real sequence.  Then for any integer $s\geq 0$, the sequence
$\{g_s(k)=\frac{1}{k}\sum_{i=s+1}^{s+k}z_i, k\geq 1\}$ is monotonically nondecreasing (nonincreasing) for $k$.
\end{lem}
\begin{lem}\label{propiid}
For i.i.d. random variables $\{\xi_i(t), i\in\mathcal{V}, t\geq 1\}$ with $E\xi_1(1)=0,\,E\xi_1^2(1)>0$, there exist constants $a>0$ and $0<p\leq 1$, such that
\begin{equation*}
  P\{\xi_i(t)\geq a\}\geq p,\,\,\,P\{\xi_i(t)\leq -a\}\geq p.
\end{equation*}
\end{lem}
Lemma \ref{propiid} is quite straightforward and we omit its proof.
\begin{lem}\label{homostublem}
Let $0<\delta\leq \frac{\epsilon}{2(n+1)}$ and suppose a.s. there is a finite time $T\geq 0$ such that $d_\mathcal{V}^{B_1}(T)\leq(n+1)\delta$, then a.s. $d_\mathcal{V}\leq2\delta, d_\mathcal{V}^{B_1}\leq (n+1)\delta$.
\end{lem}
Lemma \ref{homostublem} implies that once all the noisy opinions enter the neighbor region of stubborn agents, they will never get far away, and further reach $2\delta$-consensus. The proof of Lemma \ref{homostublem} is put in Appendix.

\noindent{\it Proof of Theorem \ref{thmhomo}:} (c) By Lemma \ref{propiid}, there exist $0<a\leq \delta, 0<p\leq 1$ such that for all $i\in \mathcal{V}, t\geq 1$,
\begin{equation}\label{probnoisevalue}
  P\{\xi_i(t)\geq a\}\geq p, \quad P\{\xi_i(t)\leq -a\}\geq p.
\end{equation}
Denote $\widetilde{x}_i(t)=|\mathcal{N}(i, x(t))|^{-1}\sum_{j\in \mathcal{N}(i, x(t))}x_j(t)$, $t\geq 0$, and this denotation remains valid for the rest of the context.
If $d_\mathcal{V}^{B_1}(0)\leq (n+1)\delta$, the conclusion holds by Lemma \ref{homostublem}. Otherwise, consider the following protocol: for all $i\in\mathcal{V}, t>0$,
\begin{equation}\label{noiseproto}
  \left\{
    \begin{array}{ll}
      \xi_i(t+1)\in[a,\delta], & \hbox{if}\quad\widetilde{x}_i(t)\leq B_1; \\
      \xi_i(t+1)\in[-\delta,-a], & \hbox{if}\quad\widetilde{x}_i(t)>B_1.
    \end{array}
  \right.
\end{equation}
By Lemma \ref{monosmlem}, (\ref{HKnoise0}), (\ref{xit0})-(\ref{neigh1}), it has under the protocol (\ref{noiseproto}) that
\begin{equation}\label{dis1}
\begin{split}
  d_\mathcal{V}^{B_1}(1)&\leq \max_{i\in\mathcal{V}}\{|\widetilde{x}_i(0)-B_1|-a\}\leq \max_{i\in\mathcal{V}}\{|x_i(0)-B_1|-a\}\\
&\leq d_\mathcal{V}^{B_1}(0)-a.
\end{split}
\end{equation}
By (\ref{probnoisevalue}) and independence, we have
\begin{equation}\label{probdistt1}
  P\{d_\mathcal{V}^{B_1}(1)\leq d_\mathcal{V}^{B_1}(0)-a\}\geq p^n.
\end{equation}
Let $L=\frac{1-(n+1)\delta}{a}$, if $d_\mathcal{V}^{B_1}(1)\leq (n+1)\delta$ a.s., the conclusion holds by Lemma \ref{monosmlem}. Otherwise, continue the above procedure $L$ times. By (\ref{probdistt1}) and independence of noises, it has
\begin{equation}\label{probdisttL}
\begin{split}
    P\{d_\mathcal{V}^{B_1}(L)\leq (n+1)\delta\}\geq &P\{\text{protocol (\ref{noiseproto}) occurs L+1 times}\}\\
    \geq &p^{n(L+1)}>0.
    \end{split}
\end{equation}
Thus
\begin{equation}\label{probdistnotL}
    P\Big\{d_\mathcal{V}^{B_1}(L)> (n+1)\delta\Big\}\leq 1-p^{n(L+1)}.
\end{equation}
Denote events $(m\geq 1)$
\begin{equation}\label{eventdelt}
\begin{split}
E_0=&\Omega, \\
E_m=&\Big\{\omega: d_\mathcal{V}^{B_1}(t)>(n+1)\delta, (m-1)L<t\leq mL\Big\}.
\end{split}
\end{equation}
Since $x(0)$ is arbitrarily given, by (\ref{probdistnotL}), it has for $m\geq 1$
\begin{equation}\label{probevent}
\begin{split}
  P\Big\{E_m\Big|\bigcap\limits_{j<m}E_j\Big\}&\leq P\Big\{d_\mathcal{V}^{B_1}(mL)> (n+1)\delta\Big|\bigcap\limits_{j<m}E_j\Big\}\\
  &\leq 1-p^{n(L+1)}<1.
  \end{split}
\end{equation}
By Lemma \ref{homostublem}, it must hold
 $\Big\{d_\mathcal{V}^{B_1}> (n+1)\delta\Big\}\subset\Big\{\bigcap\limits_{m\geq 1}\{d_\mathcal{V}^{B_1}(t)> (n+1)\delta, (m-1)L<t\leq mL\}\Big\}$, subsequently by (\ref{probevent})
\begin{equation*}
\begin{split}
P\Big\{d_\mathcal{V}^{B_1}\leq (n+1)\delta\Big\}=&1-P\Big\{d_\mathcal{V}^{B_1}> (n+1)\delta\Big\}\\
\geq &1-P\bigg\{\bigcap\limits_{m\geq 1}\{d_\mathcal{V}^{B_1}(t)> (n+1)\delta, (m-1)L<t\leq mL\}\bigg\}\\
=& 1-P\bigg\{\bigcap\limits_{m\geq 1}E_m\bigg\} =-P\bigg\{\lim\limits_{m\rightarrow \infty}\bigcap\limits_{j= 1}^mE_j\bigg\}\\
=&1-\lim\limits_{m\rightarrow \infty}P\bigg\{\bigcap\limits_{j= 1}^mE_j\bigg\} \\
=& 1-\lim\limits_{m\rightarrow \infty}\prod\limits_{j=1}^mP\bigg\{E_j\bigg|\bigcap\limits_{k<j}E_k\bigg\}\\
\geq & 1-\lim\limits_{m\rightarrow \infty}\Big(1-p^{n(L+1)}\Big)^m=1.
\end{split}
\end{equation*}
Here,  the exchangeability of probability and limit holds since $\Big\{\bigcap\limits_{j= 1}^mE_j\Big\}, m\geq 1$ is a
decreasing sequence and $P$ is a probability measure (refer to Corollary 1.5.2 of ~\cite{Chow1997}).

Since $\delta\in(0,\frac{\epsilon}{2(n+1)})$, it has $P\{d_\mathcal{V}^{B_1}< \epsilon/2\}=1$.
Consequently, a.s. there exists a finite time $T$ such that $d_\mathcal{V}^{B_1}(T)\leq \epsilon/2$. Recalling Lemma \ref{homostublem} completes the proof of (c).\hfill $\Box$

\section{HK model with heterogenous prejudices}\label{Sechetermodel}
In this part, we will mainly examine the HK model with heterogeneously prejudiced agents, while the one with heterogeneously stubborn agents will be interpreted simply at last.``Heterogeneous" means that the prejudiced agents in the group possess pronouncedly distinct opinions. Heterogeneous prejudices (or preferences) exist ubiquitously in real society, and taking rightwing and leftwing in politics as example, produce much social cleavage. In this part, we introduce the noisy HK models with heterogeneously prejudiced agents, afterwards establish the theoretical results of fragmentation phenomenon. The heterogeneous model differs from the homogeneous one in that more than one prejudiced opinion values are assumed in the models. In the subsequently proposed model, only two distinctly heterogeneous opinions are discussed for convenience.

\subsection{Noisy HK model with heterogeneously prejudiced agents}
Similarly, the heterogeneously prejudiced model is obtained by extending (\ref{xit1}):
\begin{equation}\label{xit2}
\begin{split}
  x_i^*(t)=(1-\alpha I_{\{i\in\mathcal{S}_1\cup\mathcal{S}_2\}})\sum\limits_{j\in\mathcal{N}_i(x(t))}\frac{x_j(t)}{|\mathcal{N}_i(x(t))|}+\alpha (I_{\{i\in\mathcal{S}_1\}}J_1+I_{\{i\in\mathcal{S}_2\}}J_2)+\xi_i(t+1),
\end{split}
\end{equation}
where $\mathcal{S}_1\cup\mathcal{S}_2=\mathcal{V},\mathcal{S}_1\cap\mathcal{S}_2=\emptyset$ are the sets of heterogeneously prejudiced agents, and $J_1,J_2\in[0,1]$ are the distinct prejudiced opinion values satisfying $|J_1-J_2|>\epsilon$.

\subsection{Robust fragmentation of heterogeneous prejudices}
It is easy to examine that without noise, system (\ref{HKnoise0}), (\ref{neigh0}), (\ref{xit2}) will produce fragmentation. The following theorems tell that the fragmentation phenomenon is preserved and even refined in the presence of noise.
First we present a general but conservative conclusion, then derive a more refined result.
\begin{thm}\label{thmheter}
Consider system (\ref{HKnoise0}), (\ref{neigh0}), (\ref{xit2}), given any $x(0)\in[0,1]^n, \epsilon\in(0,1)$, then (i) $\mathcal{S}_1$ a.s. achieves $\frac{(1-\alpha)\epsilon+\delta}{\alpha}$-consensus with $J_1$; (ii) $\mathcal{S}_2$ a.s. achieves $\frac{(1-\alpha)\epsilon+\delta}{\alpha}$-consensus with $J_2$.
\end{thm}
\noindent\emph{Proof}:
Consider $i\in\mathcal{S}_1$, then by (\ref{xit2}), it has for $t\geq 0$
\begin{equation*}
  x_i(t+1)=\alpha J_1+(1-\alpha)\sum\limits_{j\in\mathcal{N}_i(x(t))}\frac{x_j(t)}{|\mathcal{N}_i(x(t))|}+\xi_i(t+1).
\end{equation*}
Noting that $|x_i(t)-x_j(t)|\leq \epsilon$ for all $j\in\mathcal{N}_i(x(t))$, it obtains
\begin{equation}\label{xits1}
\begin{split}
|x_i(t+1)-J_1|=&\Big|(1-\alpha)\sum\limits_{j\in\mathcal{N}_i(x(t))}\frac{x_j(t)-J_1}{|\mathcal{N}_i(x(t))|}+\xi_i(t+1)\Big|\\
\leq& (1-\alpha)\sum\limits_{j\in\mathcal{N}_i(x(t))}\frac{|x_j(t)-J_1|}{|\mathcal{N}_i(x(t))|}+|\xi_i(t+1)|\\
\leq &(1-\alpha)\sum\limits_{j\in\mathcal{N}_i(x(t))}\frac{|x_i(t)-J_1|+\epsilon}{|\mathcal{N}_i(x(t))|}+\delta\\
=&(1-\alpha)|x_i(t)-J_1|+(1-\alpha)\epsilon+\delta\\
&\cdots\\
\leq& (1-\alpha)^{t+1}|x_i(0)-J_1|+((1-\alpha)\epsilon+\delta)\\
&\cdot((1-\alpha)^t+\ldots+(1-\alpha)+1)\\
\rightarrow &\frac{(1-\alpha)\epsilon+\delta}{\alpha},\quad a.s.,\quad as \,\,t\rightarrow\infty,
  \end{split}
\end{equation}
implying $d_{\mathcal{S}_1}^{J_1}\leq \frac{(1-\alpha)\epsilon+\delta}{\alpha}$, a.s.. Similarly, we can get that $d_{\mathcal{S}_2}^{J_2}\leq \frac{(1-\alpha)\epsilon+\delta}{\alpha}$, a.s..
\hfill $\Box$

Theorem \ref{thmheter} says that the prejudiced agents will finally stay near their preferences and the system generally displays a rough robustness of fragmentation. Smaller $\epsilon$ and $\delta$, and larger $\alpha$ make a closer reach to $J_1$ and $J_2$. Though all agents in the group admit two distinct prejudices, sometimes the noise-free fragmentation forms with more than two clusters, while in the presence of noise, the group emerges with bipartite cleavage (see Figures \ref{noisefreeprebfig} and \ref{prejudbfig}). The following theorem gives a guarantee of the bipartite fragmentation when the discrepancy of $J_1$ and $J_2$ is large enough,.
\begin{thm}\label{thmbiclu}
Suppose $J_1-J_2>\epsilon+2\frac{(1-\alpha)\epsilon+\delta}{\alpha}$ in system (\ref{HKnoise0}), (\ref{neigh0}), (\ref{xit2}), then for any $x(0)\in[0,1]^n$, (i) $\mathcal{S}_1$ a.s. achieves $\frac{\delta}{\alpha}$-consensus with $J_1$; (ii) $\mathcal{S}_2$ a.s. achieves $\frac{\delta}{\alpha}$-consensus with $J_2$.
\end{thm}
\begin{lem}\label{lembiclu}
Let $d_{\mathcal{S}_1}^{J_1}(0)\leq \frac{(1-\alpha)\epsilon+\delta}{\alpha}, d_{\mathcal{S}_2}^{J_2}(0)\leq \frac{(1-\alpha)\epsilon+\delta}{\alpha}$ in system (\ref{HKnoise0}), (\ref{neigh0}), (\ref{xit2}) with $J_1-J_2>\epsilon+2\frac{(1-\alpha)\epsilon+\delta}{\alpha}$. If a.s. there exists a finite time $T<\infty$ that $d_{\mathcal{S}_1}^{J_1}(T)\leq \frac{\delta}{\alpha}$, then $d_{\mathcal{S}_1}^{J_1}\leq \frac{\delta}{\alpha}$; Similarly, if a.s. there exists a finite time $d_{\mathcal{S}_2}^{J_2}(T)\leq \frac{\delta}{\alpha}$, then $d_{\mathcal{S}_2}^{J_2}\leq \frac{\delta}{\alpha}$.
\end{lem}
Proof of Lemma \ref{lembiclu} is put in Appendix.

\noindent \emph{Proof of Theorem \ref{thmbiclu}:}
By Theorem \ref{thmheter}, there a.s. exists $T<\infty$ such that for all $t\geq T$, it has $|x_i(t)-J_1|\leq \frac{(1-\alpha)\epsilon+\delta}{\alpha}, |x_j(t)-J_2|\leq \frac{(1-\alpha)\epsilon+\delta}{\alpha}, i\in\mathcal{S}_1, j\in\mathcal{S}_2$. Since $J_1-J_2>\epsilon+2\frac{(1-\alpha)\epsilon+\delta}{\alpha}$, we know that from the moment $T$, any agent in $\mathcal{S}_1$ cannot be the neighbor of agents in $\mathcal{S}_2$, and vice versa. Without loss of generality, we suppose $T=0$ a.s. Consider the subgroup $\mathcal{S}_1$, if $d_{\mathcal{S}_1}^{J_1}(0)\leq \frac{\delta}{\alpha}$, by Lemma \ref{lembiclu} we have $d_{\mathcal{S}_1}^{J_1}\leq \frac{\delta}{\alpha}$.
Otherwise, consider the following protocol: for all $i\in\mathcal{S}_1, t>0$,
\begin{equation}\label{noiseprotobij}
  \left\{
    \begin{array}{ll}
      \xi_i(t+1)\in[a,\delta], & \hbox{if}\quad\widetilde{x}_i(t)\leq J_1; \\
      \xi_i(t+1)\in[-\delta,-a], & \hbox{if}\quad\widetilde{x}_i(t)>J_1,
    \end{array}
  \right.
\end{equation}
where $a>0$ is given in Lemma \ref{propiid}.
By Lemma \ref{monosmlem} and (\ref{HKnoise0}), (\ref{neigh0}), (\ref{xit2}), it has under the protocol (\ref{noiseprotobij}) that
\begin{equation}\label{dis1bij}
\begin{split}
  d_{\mathcal{S}_1}^{J_1}(1)&\leq \max_{i\in\mathcal{S}_1}\{|\widetilde{x}_i(0)-J_1|-a\}\leq \max_{i\in\mathcal{S}_1}\{|x_i(0)-J_1|-a\}\\
&\leq d_{\mathcal{S}_1}^{J_1}(0)-a.
\end{split}
\end{equation}
By Lemma \ref{propiid} and independence, we have
\begin{equation}\label{probdistt1bij}
  P\{d_{\mathcal{S}_1}^{J_1}(1)\leq d_{\mathcal{S}_1}^{J_1}(0)-a\}\geq p^n.
\end{equation}
Then following a similar line of the procedure behind (\ref{probdistnotL}) of the proof of Theorem \ref{thmhomo}, we can get that
\begin{equation*}
  P\Big\{d_{\mathcal{S}_1}^{J_1}\leq \frac{\delta}{\alpha}\Big\}=1.
\end{equation*}
Similarly, we can also obtain
\begin{equation*}
  P\Big\{d_{\mathcal{S}_2}^{J_2}\leq \frac{\delta}{\alpha}\Big\}=1.
\end{equation*}
This completes the proof. \hfill $\Box$

At the end of this section, we give a short discussion of the noisy HK model with heterogeneously stubborn agents.
Introducing in (\ref{stubbornvalues1}) another stubborn agent set, $\mathcal{B}_2$,  and the distinct stubborn opinion values satisfy
\begin{equation}\label{stubbornvalues2}
  x_i(t)\equiv B_1,\quad x_j(t)\equiv B_2,\quad i\in\mathcal{B}_1, j\in\mathcal{B}_2, t\geq 0,
\end{equation}
with $B_1, B_2\in[0,1], \mathcal{B}_1\cap\mathcal{B}_2=\emptyset$. The neighbor set in (\ref{neigh1}) is modified accordingly as
\begin{equation}\label{neigh2}
  \mathcal{N}(i,x(t))=\{j\in \mathcal{V}\cup\mathcal{B}_1\cup\mathcal{B}_2:|x_j(t)-x_i(t)|\leq \epsilon\}
\end{equation}
Further, it assumes that $B_2-B_1>\epsilon$ which implies a pronounced difference between the stubborn agents.

Theorem \ref{thmbiclu} clearly shows that given any initial opinion values, the system with heterogeneous prejudices will display robust fragmentation in the presence of noise. However, for the system (\ref{HKnoise0}), (\ref{xit0}), (\ref{stubbornvalues2}), (\ref{neigh2}) of heterogeneously stubborn agents, whether robust fragmentation can be formed highly depends on the location of initial opinions. Here, without further interpretation, we simply give some examples:
\begin{thm}
For system (\ref{HKnoise0}), (\ref{xit0}), (\ref{stubbornvalues2}), (\ref{neigh2}), (i) if $x_i(0)\in[0,B_2-\epsilon)$ or $(B_1+\epsilon,B_2]$, then for all $\delta\in(0,\frac{B_2-B_1-\epsilon}{n+1})$, we have $d_\mathcal{V}\leq 2\delta$; (ii) if $x_i(0)\in[0,B_1], i\in\mathcal{V}_1$ and $x_j(0)\in[B_2,1], j\in\mathcal{V}_2$ where $\mathcal{V}_1\cup\mathcal{V}_2=\mathcal{V}, \mathcal{V}_1\cap\mathcal{V}_2=\emptyset$, then for all $\delta\in(0,\frac{B_2-B_1-\epsilon}{2(n+1)})$, we have $d_{\mathcal{V}_1}\leq 2\delta, d_{\mathcal{V}_2}\leq 2\delta$.
\end{thm}

\section{Simulations}\label{Secsimul}
\begin{figure}[htp]
  \centering
  \includegraphics[width=4.5in]{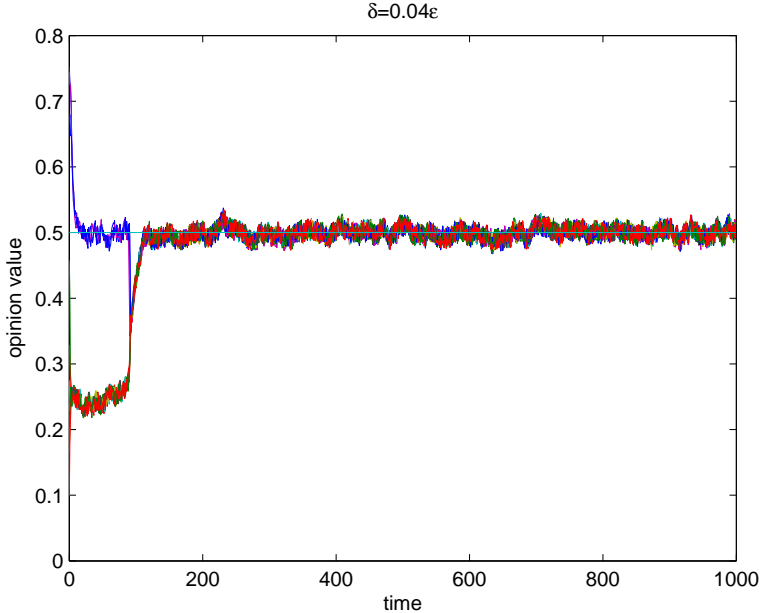}\\
  \caption{Opinion evolution of system (\ref{HKnoise0}), (\ref{xit0})-(\ref{neigh1}) of 10 agents and a stubborn agent. The initial opinion value is randomly generated
  in $[0,1]$, the opinion value of stubborn agents takes 0.5,  confidence threshold $\epsilon=0.2$, and noise strength $\delta=0.04\epsilon$.  }\label{stubbornsfig}
\end{figure}

In this part, we present some simulation results to verify our main theoretical conclusions. First consider the HK model with homogeneously stubborn agents which consists of 10 agents whose initial opinion values are randomly generated in $[0,1]$, and a stubborn agent whose constant opinion value takes 0.5. Further, let $\epsilon=0.2, \delta=0.04\epsilon$, then Figure \ref{stubbornsfig} shows that the opinions finally synchronize to 0.5 and the fragmentation vanishes. Some simulation results of the systems with homogeneously prejudiced agents can be found in \cite{Su2017free}, and we omit them here.

\begin{figure}[htp]
  \centering
  \includegraphics[width=4.5in]{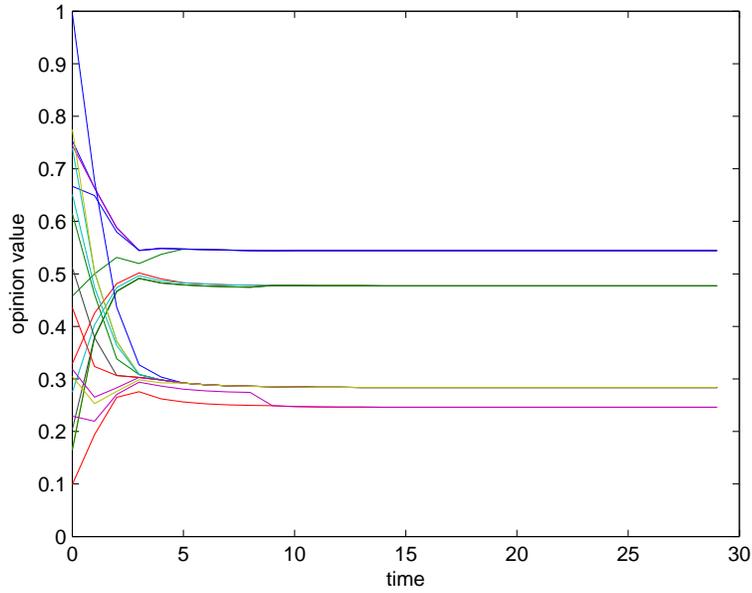}\\
  \caption{Opinion evolution of system (\ref{HKnoise0}), (\ref{neigh0}), (\ref{xit2}) of 20 agents without noise. The initial opinion value is randomly generated in $[0,1]$
  ,$J_1=0.6, J_2=0.2$, confidence threshold $\epsilon=0.2$, and attraction strength $\alpha=0.4$.  }\label{noisefreeprebfig}
\end{figure}

\begin{figure}[htp]
  \centering
  \includegraphics[width=4.5in]{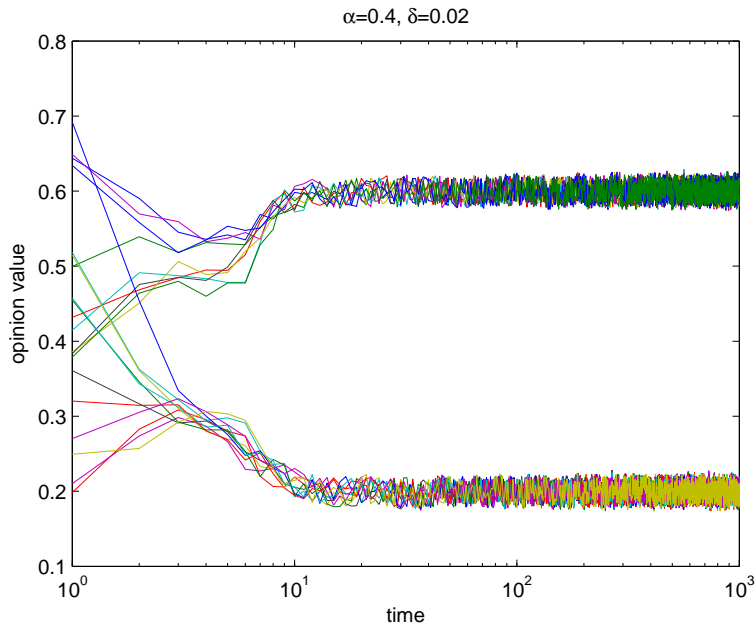}\\
  \caption{Opinion evolution of system (\ref{HKnoise0}), (\ref{neigh0}), (\ref{xit2}) of 20 agents with noise. The initial conditions are the same as that in Figure \ref{noisefreeprebfig} and the noise strength $\delta=0.02$.  }\label{prejudbfig}
\end{figure}

In the following, we consider the HK system with heterogeneously prejudiced agents. Let the system consist of 20 agents, whose initial opinion values randomly generate in $[0,1]$, $J_1=0.6, J_2=0.2$, confidence threshold $\epsilon=0.2$, and attraction strength $\alpha=0.4$. First we present the opinion evolution without noise, and Figure \ref{noisefreeprebfig} shows that the system forms four clusters. Then we add noise with strength $\delta=0.02$ in the system, and Figure \ref{prejudbfig} shows that the fragmentation emerges with two clusters which locate near the prejudices $J_1$ and $J_2$.

\section{Conclusions}\label{Secconc}
In this paper, we investigated how the fragmentation phenomenon emerges with HK dynamics. It has been shown that the original HK model fails to generate cleavage in the presence of tiny noise. Here we further prove that the well-appearing fragmentation of HK models with homogeneously prejudiced or stubborn agents also vanishes with noise. We finally reveal that the HK model with heterogeneously prejudiced agents can preserve a robust and well-shaping fragmentation under the drive of noise. This implies a quite intuitive conclusion that only the innate difference within the group rather than the mere bounded confidence of individuals that cause the ubiquitous social cleavage.

\vspace{1.5ex}
\noindent\textbf{Appendix}

\vspace{1.5ex}
\noindent{\it Proof of Lemma \ref{homostublem}:} It is easy to check that $d_\mathcal{V}(T)\leq 2|d_\mathcal{V}^{B_1}(T)|\leq \epsilon$. This implies that at time $T$, all agents in $\mathcal{V}$ are neighbors to each other. By (\ref{HKnoise0}), (\ref{xit0})-(\ref{neigh1}), we have for $i\in\mathcal{V}$
\begin{equation*}
  x_i(T+1)= \frac{\sum\limits_{j\in\mathcal{V}}x_j(T)+|\mathcal{B}_1|B_1}{n+|\mathcal{B}_1|}+\xi_i(T+1)
\end{equation*}
thus
\begin{equation}\label{neihomolem}
  \begin{split}
    |x_i(T+1)-B_1|\leq & \bigg|\frac{\sum_{j\in\mathcal{V}}x_j(T)+|\mathcal{B}_1|B_1}{n+|\mathcal{B}_1|}-B_1\bigg|+|\xi_i(T+1)| \\
      \leq & \frac{\sum_{j\in\mathcal{V}}|x_j(T)-B_1|}{n+|\mathcal{B}_1|}+\delta\leq \frac{n(n+1)}{n+|\mathcal{B}_1|}\delta+\delta\\
\leq &(n+1)\delta\leq \frac{\epsilon}{2},\,\,a.s..
  \end{split}
\end{equation}
Meanwhile, for all $i,j\in\mathcal{V}$, it has a.s.
\begin{equation*}
\begin{split}
  |x_i(T+1)-x_j(T+1)|\leq |\xi_i(T+1)|+|\xi_j(T+1)|\leq 2\delta.
  \end{split}
\end{equation*}
Repeating the above procedure yields the conclusion. \hfill $\Box$

\noindent\emph{Proof of Lemma \ref{lembiclu}:}
Since $J_1-J_2>\epsilon+2\frac{(1-\alpha)\epsilon+\delta}{\alpha}$, at time $T$, any agent in $\mathcal{S}_1$ cannot be the neighbor of agents in $\mathcal{S}_2$, hence for all $i\in\mathcal{S}_1$
\begin{equation*}
\begin{split}
  |x_i(T+1)-J_1|=&\Big|(1-\alpha)\sum\limits_{j\in\mathcal{N}_i(x(T))}\frac{x_j(T)-J_1}{|\mathcal{N}_i(x(T))|}+\xi_i(T+1)\Big|\\
\leq& (1-\alpha)\sum\limits_{j\in\mathcal{N}_i(x(T))}\frac{|x_j(T)-J_1|}{|\mathcal{N}_i(x(T))|}+|\xi_i(T+1)|\\
\leq &(1-\alpha)\frac{\delta}{\alpha}+\delta=\frac{\delta}{\alpha},\,\,a.s.
  \end{split}
\end{equation*}
Repeating the above procedure, it follows that $d_{\mathcal{S}_1}^{J_1}\leq \frac{\delta}{\alpha}$, and similarly $d_{\mathcal{S}_2}^{J_2}\leq \frac{\delta}{\alpha}$.
\hfill $\Box$

\end{document}